\documentclass[11pt]{article}
\usepackage{amssymb}
\usepackage{mathrsfs}
\usepackage{amsfonts}
\usepackage{tikz}

\newtheorem {Problem} {Problem}[section]
\newtheorem {Theorem} [Problem]{Theorem}
\newtheorem {Lemma}[Problem]{Lemma}

\newtheorem {Corollary}[Problem]{Corollary}
\newenvironment {Proof}{\noindent {\bf Proof.}}{\hfill\ensuremath{\square}}
\newcommand*{\QEDB}{\hfill\ensuremath{\square}}

\begin{document}

\title{ Spectral extremal results with forbidding linear forests \thanks{This work is supported by  the Joint NSFC-ISF Research Program (jointly funded by the National Natural Science Foundation of China and the Israel Science Foundation (No. 11561141001)) and  the National Natural Science Foundation of China (No.11531001)
}
}

\author{ Ming-Zhu Chen, A-Ming Liu,
 Xiao-Dong Zhang\footnote{Corresponding author. E-mail: xiaodong@sjtu.edu.cn}
 \\
School of Mathematics Science, MOE-LSC, SHL-MAC\\
Shanghai Jiao Tong University,
Shanghai 200240, P. R. China}

\date{}
\maketitle

\begin{abstract}
 The Tur\'{a}n type extremal problem  asks to maximize the number of edges over all graphs which do not
contain fixed  subgraphs. Similarly, the spectral Tur\'{a}n type extremal problem  asks to maximize spectral radius of all graphs which do not contain fixed subgraphs.
 In this paper, we  determine the maximum spectral radius of  all  graphs  without  containing a linear  forest as a subgraph  and  characterize all corresponding extremal graphs. In addition, the maximum number of edges and spectral radius of all bipartite graphs without containing $k\cdot P_3$  as a subgraph  are obtained and all extremal graphs are also characterized. Moreover, some relations between Tu\'{a}n type extremal problems and spectral  Tur\'{a}n type extremal problems are discussed.
\end{abstract}

\section{Introduction}
 Let $G$ be an undirected simple graph with vertex set
$V(G)=\{v_1,\dots,v_n\}$ and edge set $E(G)$, where  $e(G)$ is the number of
edges of $G$.
The \emph{adjacency matrix}
$A(G)=(a_{ij})$ of $G$  is the $n\times n$ matrix, where
$a_{ij}=1$ if $v_i$ is adjacent to $v_j$, and $0$ otherwise.   The \emph{spectral radius} of $G$ is the largest eigenvalue of $A(G)$, denoted by
$\rho(G)$, while  the least eigenvalue of $A(G)$ is denoted by $\lambda_n(G)$.
 For $v\in V(G)$,  the \emph{neighborhood} $N_G(v)$ of $v$  is $\{u: uv\in E(G)\}$ and the \emph{degree} $d_G(v)$ of $v$  is $|N_G(v)|$.
We write $N(v)$ and $d(v)$ for $N_G(v)$ and $d_G(v)$ respectively if there is no ambiguity. A path of order $n$ is denoted by $P_n$.
 For $V_1, V_2\subseteq V(G)$,  $e(V_1,V_2)$ denotes the number of the  edges of  $G$ with one end vertex in $V_1$ and the other  in $V_2$.
A graph $G$ is said \emph{$F$-free} if   it  does not contain $F$ as a subgraph.  A {\it linear forest} is a forest whose connected components are paths. For a path $P_3$ of order $3$, say $xyz$, we call $y$ its center and $x,z$ its two ends.
 For two vertex disjoint graphs $G$ and $H$,  we denote by  $G\cup H$ and  $G\vee H$  the \emph{union} of $G$ and $H$, and the \emph{join} of $G$ and $H$, i.e., joining every vertex of $G$ to every vertex of $H$, respectively.
Denote by $k\cdot G$  the $k$ disjoint union of $G$.
 For graph notation and terminology undefined here, we refer the readers to \cite{BM}.

The problem of maximizing the number of edges over all graphs without containing fixed subgraphs is one of the cornerstones of graph theory.
In 2010,  Nikiforov \cite{Nikiforov} proposed the following
 spectral extremal graph  problem,  which is the spectral analogue of Tur\'{a}n type extremal problem.

\begin{Problem}\label{P1}
Given a graph $H$, what is the maximum $\rho(G)$ of a graph $G$ of order $n$ which does not contain $H$ as a subgraph$?$
\end{Problem}

 On Problem \ref{P1}, Nikiforov has obtained a bulk of work in this spectral analogue of Tur\'{a}n type extremal problem.
 For example, he presented some  spectral analogues of classical  results in extremal graph theory, such as  spectral analogue of  Tur\'{a}n type theorem \cite{Nikiforov2} and the Erd\H{o}s--Stone--Bollob\'{a}s theorem \cite{Nikiforov3}.
For more details,  readers may be referred to \cite{Nikiforov4,Nikiforov5,Nikiforov6,Nikiforov,TT,YWZ,ZW}.
In particular, Nikiforov \cite{Nikiforov} determined the maximum spectral radius of a graph that does not contain  paths
of given length as  subgraphs and characterized all  extremal graphs, which is a spectral analogue of Erd\H{o}s-Gallai theorem (see \cite{EG}).

 A natural extension of this problem is to determine the maximum spectral radius of graphs without containing a linear forest.
 Lidick\'{y}, Liu, and Palmer \cite{LLP}  determined the Tur\'{a}n number for a forbidden linear forest except for $ k\cdot P_3$  if the order of graph is  sufficiently large  and characterized all extremal graphs.  Bushaw and Kettle \cite{BK}, Campos and Lopes \cite{campos2015}, and Yuan and Zhang \cite{YZ}, independently,  determined  the Tur\'{a}n number for a forbidden $ k\cdot P_3$. In order to state these results, we need some symbols for given graphs.

 For $0\leq s\leq \lfloor\frac{n-1}{2}\rfloor$,  let $T_{n,s}$ be the graph of order $n$ obtained by identifying an end  of each of $s$ paths $P_3$ and a vertex of each of $n-2s-1$ paths $P_2$. Clearly, $T_{n,0}$ is a star of order $n$.

For $0<h<n, $ let $S_{n,h}$ be the graph of order $n$ obtained by joining every vertex of  the complete graph $K_h$ of order $h$  to every vertex of  the complement graph  $\overline{K}_{n-h}$ of $K_{n-h}$, i.e., $S_{n,h}=K_h\vee \overline{K}_{n-h}$.
Furthermore, let $S^+_{n,h}$ be the graph obtained by adding an edge to $S_{n,h}$, i.e., $S^+_{n,h}=K_2\vee(K_h\cup K_{n-h-2})$.

For $1\le k<n$, let $F_{n,k}:=K_{k-1}\vee ( pK_2\cup K_s)$, where $n-(k-1)=2p+s$ and $0\leq s<2$. In particular,
$F_{n, 1}= pK_2\cup K_s$.

\begin{Theorem}\label{edge extremal graph for linear forest}\cite{LLP}
Let  $F$ be a linear forest, i.e., $F=\cup_{i=1}^k P_{a_i}$,  with $k\ge 2 $, $a_1\geq \cdots \geq a_k\geq2$,  and  $h=\sum \limits_{i=1}^k \lfloor\frac{a_i}{2}\rfloor-1$. If   there exists at least one $a_i$ not 3 and $G$ is  an $F$-free  graph  of  order $n$,  then for  sufficiently large $n$,
$$e(G)\leq   \left(\begin{array}{c} h\\ 2\end{array}\right)  + h(n-h)+c,$$
where $c=1$ if all $a_i$ are odd and $c=0$ otherwise. Moreover,
If $c=1$ then  equality holds if and only if $G=S^+_{n,h}$. Otherwise, the equality holds if and only if  $G=S_{n,h}$.
\end{Theorem}


\begin{Theorem}\label{edge extremal graph for kP3}\cite{BK, campos2015,YZ}
Let  $G$ be a \ $k\cdot P_3$-free graph of order $n$. Then 
\[e(G)\leq\left\{
\begin{array}{llll}
 \vspace{1mm}
 \left(\begin{array}{c} n\\ 2\end{array}\right) ,&& \mbox{for $n<3k$};\\
   \vspace{1mm}
 \left(\begin{array}{c} 3k-1\\ 2\end{array}\right) +\big\lfloor\frac{n-3k+1}{2}\big\rfloor,&& \mbox{for $3k\leq n<5k-1$};\\
       \vspace{1mm}
  \left(\begin{array}{c} 3k-1\\ 2\end{array}\right) +k, && \mbox{for $n=5k-1$};\\
      \vspace{1mm}
 \left(\begin{array}{c} k-1\\ 2\end{array}\right)+(n-k+1)(k-1)+\big\lfloor\frac{n-k+1}{2}\big\rfloor, && \mbox{for $n>5k-1$.}
\end{array}\right.
\]
Moreover, (i) If $n <3k$, then equality holds if and only if $G=K_n$;\\
(ii) If $3k\leq n < 5k-1$, then  equality holds if and only if  $G=K_{3k-1}\cup F_{n-3k+1,1}$;\\
(iii) If $n = 5k - 1$, then equality holds if and only if  $G=K_{3k-1}\cup F_{2k,1}$ or $G=F_{5k-1,k}$;\\
(iv) If $n > 5k - 1$, then equality holds if and only if $G=F_{n,k}$.
\end{Theorem}

Motivated by Problem \ref{P1} and  above results,  we determine the  maximum spectral  radius of all graphs  without a  linear forest and characterize all extremal graphs. In addition, we also determine the maximum number of edges and spectral radius of  bipartite graphs which do not contain $k\cdot P_3$  as a subgraph
  and   characterize all extremal graphs. The main results of this paper are stated as follows.

  \begin{Theorem}\label{spectral radius for all linear forest}
 Let  $F$ be a linear forest, i.e., $F=\cup_{i=1}^k P_{a_i}$  with   $k\geq2$ and
$a_1\geq \cdots \geq a_k\geq2$. Denote $ h=\sum \limits_{i=1}^k \lfloor\frac{a_i}{2}\rfloor-1$ and suppose that $G$ is  an $F$-free  graph  of sufficiently large order $n$. \\
(i) If there  exists an even $a_i$, then $\rho(G)\leq \rho(S_{n,h})$ with equality if and only if $G=S_{n,h}$;\\
(ii) If all $a_i$ are odd  and there exists at least one $a_i>3$,   then $\rho(G)\leq \rho(S^+_{n,h})$ with equality if and only if $G=S^+_{n,h}$.\\
(iii) If all $a_i$ are 3, i.e., $F=k\cdot P_3$, then $\rho(G)\leq \rho(F_{n,k}) $ with equality if and
only if   $G=F_{n,k}$.
\end{Theorem}

\begin{Theorem}\label{edge extremal bipartite graph}
Let  $G$ be a $k\cdot P_3$-free bipartite graph of order $n\geq 11k-4$ with $k\geq2. $ Then $$e(G)\leq (k-1)(n-k+1).$$
Moreover, if  $k=2$ then equality holds  if and only if
$G=T_{n,s}$, $s=0,1,\dots,\lfloor\frac{n-1}{2}\rfloor$; if $k\geq3$ then
 equality holds if and only if $G=K_{k-1,n-k+1}$.
\end{Theorem}

\begin{Theorem}\label{extremal bipartite graph for kP3}
Let  $G$ be a $k\cdot P_3$-free bipartite graph of order $n\geq 11k-4$ with $k\geq2 $. Then $$\rho(G)\leq \sqrt{(k-1)(n-k+1)}$$ with equality if and
only if   $G=K_{k-1,n-k+1}$.
\end{Theorem}

\begin{Corollary}\label{least eigenvalue}
Let $G$ be a $k\cdot P_3$-free  graph of order $n\geq 11k-4$ with $k\geq2 $. Then $$\lambda_n(G)\geq -\sqrt{(k-1)(n-k+1)}$$
 with equality if and only if   $G=K_{k-1,n-k+1}$.
\end{Corollary}
 The rest of this paper is organized as follows. In Section 2, some known results and lemmas are presented. In Sections 3 and 4, we present the proof of Theorems~\ref{spectral radius for all linear forest}, \ref{edge extremal bipartite graph}, and Theorem~\ref{extremal bipartite graph for kP3} and Corollary~\ref{spectral radius for kP3}, respectively.  In Section 5, some relations between the Tur\'{a}n theorem and  spectral Tur\'{a}n theorem are discussed.
\section{Preliminary}
In this section, we present some known results.
\begin{Theorem}\label{upper bound}\cite{HSF, Nikiforov1}
Let $G$ be a graph of order $n$ with the minimum degree $\delta=\delta(G)$ and $e=e(G)$.
Then
$$\rho(G)\leq\frac{\delta-1+\sqrt{8e-4\delta n+(\delta+1)^2}}{2}.$$
\end{Theorem}
It is easy to see that for $2e\leq  n(n -1)$,  the function $$f(x)=\frac{x-1+\sqrt{8e-4x n+(x+1)^2}}{2}$$
is decreasing with respect to $x$ \cite{Nikiforov1}.

\begin{Theorem}\label{upper bound for bipartite graph}\cite{BFP}
Let $G$ be a bipartite graph. Then \\
$$\rho(G)\leq \sqrt{e(G)},$$
with equality if and only if $G$ is a disjoint union of a complete bipartite graph and isolated vertices.
\end{Theorem}

\begin{Lemma}\label{upper bound for S+nk}\cite{Nikiforov}
(i) For $h\geq1$ and $n>h$, $\rho(S_{n,h})=\frac{h-1+\sqrt{4hn-(3h^2+2h-1)}}{2}$.\\
(ii) For $h\geq2$ and $n\geq 4^h$, $\rho(S^+_{n,h})<\frac{h-1+\sqrt{4hn-(3h^2+2h-3)}}{2}$.
\end{Lemma}

The following lemma is a little different from its original form\cite[Lemma 14]{Nikiforov}, but it can be seen from its original proof.
\begin{Lemma}\label{minimum degree less than k}\cite{Nikiforov}
Let $c\geq0$, $h\geq2$, $n\geq 2^{4h}$, and let $G$ be a graph of order $n$. If $\delta(G)<h$
and  $$\rho(G)\geq \frac{h-1+\sqrt{4hn-4h^2+c}}{2},$$
then there exists a  subgraph $H$ of order $p\geq \lfloor\sqrt{n}\rfloor$ satisfying one of the following conditions:\\
(i) $p=\lfloor\sqrt{n}\rfloor$ and $\rho(H)>\sqrt{(2h+1)p}$;\\
(ii) $p>\sqrt{n}$, $ \delta(H)\geq h$ and
$\rho(H)>\frac{h-1+\sqrt{4hp-4h^2+c+2}}{2}$.
\end{Lemma}

\begin{Lemma}\label{upper bound for Fnk}
(i)If $n-k+1$ is even, then $\rho(F_{n,k})= \frac{k-1+\sqrt{4(k-1)n-(3k^2-2k-5)}}{2}$.\\
(ii) If $n-k+1$ is odd, then $\rho$ is the largest root of $x^3-(k-1)x^2-[(k-1)n-(k^2-k-1)]x+k-1=0$.\\
(iii)  $\frac{k-1+\sqrt{4(k-1)n-(3k^2-2k-1)}}{2}<\rho(F_{n,k})\leq\frac{k-1+\sqrt{4(k-1)n-(3k^2-2k-5)}}{2}$.
\end{Lemma}

\begin{Proof}
Denote $\rho=\rho(F_{n,k})$. Let $\mathbf x$ be a positive eigenvector of $A(F_{n,k})$ corresponding to $\rho$.
Let $n-(k-1)=2p+s$ with $0\le s<2$. By symmetry, all vertices of  subgraphs $K_{k-1}$, $pK_2$, or $K_s$ in  $F_{n,k}=K_{k-1}\vee (pK_2 \cup K_s)$ have the same  eigenvector components respectively, which are denoted by $x_1$, $x_2$, $x_3$, respectively.

(i) If $n-k+1$ is even, then $s=0$.
By  $A\mathbf x=\rho\mathbf x$, it is easy to see that
\begin{eqnarray*}
    \rho x_1 &=& (k-2)x_1 +(n-k+1)x_2, \\
    \rho x_2 &=& (k-1)x_1+x_2.
\end{eqnarray*}
 It is easy to see that
 $$\rho= \frac{k-1+\sqrt{4(k-1)n-(3k^2-2k-5)}}{2}.$$

 (ii)  If $n-k+1$ is odd, then $s=1$. By  $A\mathbf x=\rho\mathbf x$, it is easy to see that
\begin{eqnarray*}
    \rho x_1 &=& (k-2)x_1 +(n-k)x_2+x_3, \\
    \rho x_2 &=& (k-1)x_1+x_2,    \\
     \rho x_3 &=& (k-1)x_1.
\end{eqnarray*}
 Hence we have
    $$\rho^3-(k-1)\rho^2-[(k-1)n-(k^2-k-1)]\rho+k-1=0.$$
 Therefore $\rho$ is the largest root of $x^3-(k-1)x^2-[(k-1)n-(k^2-k-1)]x+k-1=0$.
In addition,   $$\rho^2-(k-1)\rho-[(k-1)n-(k^2-k-1)]=-\frac{k-1}{\rho}<0,$$
 which implies that
 $$\rho<\frac{k-1+\sqrt{4(k-1)n-(3k^2-2k-5)}}{2}.$$
Moreover, noting that $\rho> \rho(K_k)=k-1$, we have
$$\rho^2-(k-1)\rho-[(k-1)n-(k^2-k)]=\frac{\rho-(k-1)}{\rho}>0,$$
which implies that $$\rho> \frac{k-1+\sqrt{4(k-1)n-(3k^2-2k-1)}}{2}.$$
\end{Proof}

\section{Proof of Theorem~\ref{spectral radius for all linear forest}}
In order to prove Theorem~\ref{spectral radius for all linear forest}, we first prove the following three  lemmas.

  \begin{Lemma}\label{spectral radius for linear forest1}
 Let  $F$ be a linear forest, i.e., $F=\cup_{i=1}^k P_{a_i}$  with   $k\geq2$ and
$a_1\geq \cdots \geq a_k\geq2$. Denote $ h=\sum \limits_{i=1}^k \lfloor\frac{a_i}{2}\rfloor-1$ and suppose that $G$ is  an $F$-free  graph  of sufficiently large order $n$.  If there  exists an even $a_i$, then $\rho(G)\leq \rho(S_{n,h})$ with equality if and only if $G=S_{n,h}$.
\end{Lemma}
\begin{Proof}
 Let $G$ be an $F$-free graph of order $n$ with the maximum spectral radius. Set $\delta=\delta(G)$ and $e=e(G)$.
Since $S_{n,h}$ is $F$-free, by Lemma~\ref{upper bound for S+nk} (i) we have $$\rho(G)\geq \rho(S_{n,h})=\frac{h-1+\sqrt{4hn-(3h^2+2h-1)}}{2}.$$

First we assume that $h=1$. Then $F=2\cdot P_2$ or $F=P_2\cup P_3$. Obviously $$\rho(G)\geq \sqrt{n-1}.$$
By \cite[Theorem~2]{Nikiforov}, $G$ contains a $P_4$ unless $G=S_{n,1}$. Note that $G$ is $F$-free. If $F= 2\cdot P_2$, then $G=S_{n,1}$.  If  $F=P_2\cup P_3$,
then $G$ must be $P_4\cup \overline{K}_{n-4}$ or $S_{n,1}$.  However, $$\rho(P_4\cup \overline{K}_{n-4})<2\leq\sqrt{n-1}.$$
Thus  $G=S_{n,1}$.

So we now assume that $h\geq2$ and  consider the following two cases.

\vspace{2mm}
{\bf Case 1:} $\delta\geq h$. By Theorem~\ref{upper bound},
 $$ \rho(G) \leq \frac{\delta-1+\sqrt{8e-4\delta n+(\delta+1)^2}}{2} \leq \frac{h-1+\sqrt{8e-4h n+(h+1)^2}}{2}.$$
Thus $e\geq hn-\frac{h^2+h}{2}=e(S_{n,h})$. Note that $G$ is $F$-free. By Theorem~\ref{edge extremal graph for linear forest}, $G=S_{n,h}$.

\vspace{2mm}
{\bf Case 2:} $\delta<h$.
Note that $$\rho(G)\geq \frac{h-1+\sqrt{4hn-(3h^2+2h-1)}}{2}=\frac{h-1+\sqrt{4hn-4h^2+(h-1)^2}}{2}.$$
 By Lemma~\ref{minimum degree less than k},  there exists a  graph $H$ of order $p$ such that
either $p=\lfloor\sqrt{n}\rfloor$ and  $\rho(H)>\sqrt{(2h+1)p};$ or  $p>\sqrt{n}$, $ \delta(H)\geq h$ and  $$\rho(H)>\frac{h-1+\sqrt{4hp-4h^2+(h-1)^2}}{2}.$$
If $p=\lfloor\sqrt{n}\rfloor$ and  $\rho(H)> \sqrt{(2h+1)p}$, then
  $$2e(H)=tr(A^2(H))\geq\rho^2(H)> (2h+1)p>2e(S_{p,h}).$$
By Theorem~\ref{edge extremal graph for linear forest}, $G$ contains $F$ as a subgraph,  which is a contradiction.
So we now assume that $p>\sqrt{n}$, $\delta(H)\geq h$ and $$\rho(H)>\frac{h-1+\sqrt{4hp-4h^2+(h-1)^2}}{2}.$$
Applying Theorem~\ref{upper bound} again, we have
$$ \rho(H) \leq \frac{h-1+\sqrt{8e(H)-4h p+(h+1)^2}}{2}.$$
Hence $e(H)>hp-\frac{h^2+h}{2}=e(S_{p,h})$.   By Theorem~\ref{edge extremal graph for linear forest}, $H$ contains $F$ as a subgraph, which is a contradiction. So we finish the proof.
\end{Proof}

  \begin{Lemma}\label{spectral radius for all linear forest2}
 Let  $F$ be a linear forest, i.e., $F=\cup_{i=1}^k P_{a_i}$  with   $k\geq2$ and
$a_1\geq \cdots \geq a_k\geq2$. Denote $ h=\sum \limits_{i=1}^k \lfloor\frac{a_i}{2}\rfloor-1$ and suppose that $G$ is  an $F$-free  graph  of sufficiently large order $n$.
 If all $a_i$ are odd  and there exists at least one $a_i>3$,   then $\rho(G)\leq \rho(S^+_{n,h})$ with equality if and only if $G=S^+_{n,h}$.\\
\end{Lemma}
\begin{Proof}
 Let $G$ be an $F$-free graph of order $n$ with maximum spectral radius. Set $\delta=\delta(G)$ and $e=e(G)$.
Since $S^+_{n,h}$ is $F$-free, by Lemma~\ref{upper bound for S+nk} (i) we have $$\rho(G)\geq  \rho(S^+_{n,h})> \rho(S_{n,h})=\frac{h-1+\sqrt{4hn-(3h^2+2h-1)}}{2}.$$
We consider the following two cases.

\vspace{2mm}
{\bf Case 1:}  $\delta\geq h$. By Theorem~\ref{upper bound}, we have
 $$ \rho(G) \leq \frac{\delta-1+\sqrt{8e-4\delta n+(\delta+1)^2}}{2} \leq \frac{h-1+\sqrt{8e-4h n+(h+1)^2}}{2}.$$
Hence $e\geq hn-\frac{h^2+h}{2}+1=e(S^+_{n,h})$. Note that $G$ is $F$-free.
By Theorem~\ref{edge extremal graph for linear forest}, $G=S^+_{n,h}$.

\vspace{2mm}
{\bf Case 2:}  $\delta<h$.
Note that $$\rho(G)> \frac{h-1+\sqrt{4hn-(3h^2+2h-1)}}{2}=\frac{h-1+\sqrt{4hn-4h^2+(h-1)^2}}{2}.$$
 By Lemma~\ref{minimum degree less than k},  there exists a  graph $H$ of order $p$ such that
either $p=\lfloor\sqrt{n}\rfloor$ and $\rho(H)>\sqrt{(2h+1)p};$ or  $p>\sqrt{n}$, $ \delta(H)\geq h$ and $$\rho(H)>\frac{h-1+\sqrt{4hp-4h^2+(h-1)^2+2}}{2}.$$

If  $p=\lfloor\sqrt{n}\rfloor$ and $\rho(H)> \sqrt{(2h+1)p}$, then
  $$2e(H)=tr(A^2(H))\geq\rho^2(H)> (2h+1)p>2e(S^+_{p,h}).$$
By Theorem~\ref{edge extremal graph for linear forest}, $G$ contains $F$ as a subgraph,  which is a contradiction.
Now we  assume that $p> \sqrt{n}$, $\delta(H)\geq h$ and $$\rho(H)>\frac{h-1+\sqrt{4hp-4h^2+(h-1)^2+2}}{2}=\frac{h-1+\sqrt{4hp-(3h^2+2h-3)}}{2}.$$
By Lemma~\ref{upper bound for S+nk}, $\rho(H)> \rho(S^+_{p,h})$. On the other hand, by Theorem~\ref{upper bound} and $\delta(H)\ge h$,
 we have
 \begin{eqnarray*}
&& \frac{h-1+\sqrt{4hp-4h^2+(h-1)^2+2}}{2}\\
&<& \rho(H)\\
&\le&    \frac{\delta(H)-1+\sqrt{8e(H)-4\delta(H) p+(\delta(H)+1)^2}}{2}\\
&\le &\frac{h-1+\sqrt{8e(H)-4hp+(h+1)^2}}{2}.  \end{eqnarray*}
Therefore, $e(H)>hp-\frac{h^2+h}{2}+\frac{1}{4}$. So $e(H)\ge e(S^+_{p, h})$. By Theorem~\ref{edge extremal graph for linear forest} and $H$ being $F$-free,  we have $H=S^+_{p,h}$, which contradicts to $\rho(H)>\rho(S^+_{p,h})$. This completes the proof.
\end{Proof}

\begin{Lemma}\label{spectral radius for kP3}
Let  $G$ be a $k\cdot P_3$-free graph of order $n\geq 8k^2-3k$ with $k\geq2 $. Then $$\rho(G)\leq \rho(F_{n,k}) $$ with equality if and
only if   $G=F_{n,k}$.
\end{Lemma}

\begin{Proof}
Let $G$ be a $k\cdot P_3$-free graph of order $n$ with maximum spectral radius. Since $F_{n,k}$ is  $k\cdot P_3$-free, we have $\rho(G)\geq \rho(F_{n,k})$.
We first prove the following claim.

\vspace{2mm}
{\bf Claim:} There exists a vertex $u\in V(G)$ such that $d(u)\geq k$ and
$$\sum\limits_{v\in N(u)} d(v)\geq (k-1)d(u)+(k-1)n-(k^2-k-1).$$

Suppose that $\sum_{v\in N(u)} d(v)\leq (k-1)d(u)+(k-1)n-(k^2-k)$ for all $u\in V(G)$.
Let $B=A^2(G)-(k-1)A(G)-[(k-1)n-(k^2-k)]I$
 and $B_u$ be the row sum of $B$ on $u$, where $I$ is an identity matrix.
Then
$$\rho(B)=\rho^2(G)-(k-1)\rho(G)-[(k-1)n-(k^2-k)].$$
On the other hand, by \cite[Lemma~2.1]{EZ},
 $$\rho(B)\leq \max_{u\in V(G)} B_u=\max_{u\in V(G)}\bigg\{\sum\limits_{v\in N(u)} d(v)-(k-1)d(u)-[(k-1)n-(k^2-k)]\bigg\}\leq0.$$
 Hence  $\rho(G)$ is no more than the largest root of $x^2-(k-1)x-[(k-1)n-(k^2-k)]=0$, i.e.,
$$\rho(G)\leq\frac{k-1+\sqrt{4(k-1)n-(3k^2-2k-1)}}{2},$$
which contradicts to   $\rho(G)> \frac{k-1+\sqrt{4(k-1)n-(3k^2-2k-1)}}{2}$ by  Lemma~\ref{upper bound for Fnk} (iii).
In addition, since
\begin{eqnarray*}
  0 &\leq&\sum\limits_{v\in N(u)} d(v)- (k-1)d(u)-[(k-1)n-(k^2-k-1)]\\
   &\leq& d(u)(n-1)-(k-1)d(u)-[(k-1)n-(k^2-k-1)] \\
   &=& d(u)(n-k)-(k-1)(n-k)-1 \\
   &<& (d(u)-k+1)(n-k),
\end{eqnarray*}
we have $d(u)>k-1$, i.e., $d(u)\geq k$. So the claim holds.

Next we consider the following two cases.

\vspace{2mm}

{\bf Case 1:} $k\leq d(u)\leq 3k-2$.  Let $W=V(G)\backslash (\{u\}\cup N(u))$ and $C\subseteq N(u)$ be the vertex subset such that every vertex in $C$
has at least $2k$ neighbours in $W$.  We claim that $|C|= k-1$.
Indeed, if $|C|\geq k$ then we can embed $k\cdot P_3$ with all centers in $C$ and all ends in $W$ into $G$, a contradiction. If
 $|C|\leq k-2$, then we have
\begin{eqnarray*}
 \sum\limits_{v\in N(u)} d(v) &=& \sum\limits_{v\in C} d(v)+\sum\limits_{v\in N(u)\backslash C} d(v) \\
   &\leq& |C|(n-1)+(d(u)-|C|)(2k-1+d(u)) \\
   &=& (n-d(u)-2k)|C|+(2k+d(u)-1)d(u)\\
   &\leq& (n-d(u)-2k)(k-2)+(2k+d(u)-1)d(u) \\
   &\leq& (k-1)d(u)+(k-1)n-(k^2-k-1)-[n-3k(3k-2)+k^2-3k+1]\\
   &=&(k-1)d(u)+(k-1)n-(k^2-k-1)-(n-8k^2+3k+1) \\
   & < &(k-1)d(u)+(k-1)n-(k^2-k-1),
\end{eqnarray*}
which  contradicts the claim. So $|C|=k-1$.

Since $|C|=k-1$, we can also  embed $(k-1)\cdot P_3$ with all centers in $C$ and all ends in $W$ into $G-u$. Denote by  $\bigcup_{1\leq i\leq k-1} Q_i$ the  $(k-1)\cdot P_3$ embedded  into $G-u$, where $Q_i=x_iy_iz_i$, $y_i\in C$, $x_i,z_i\in W$ for $1\leq i\leq k-1$. We claim that $d(u)=k$. Otherwise  $\bigcup _{1\leq i\leq k-1} Q_i$ together with a disjoint $P_3$ with center $u$ and two ends in $N(u)\backslash C$  will yield  $k\cdot P_3$, a contradiction. Then there exists  exactly one vertex $y_k\in N(u)\setminus C$.
Let $W_1=\{x_1,\dots,x_{k-1},z_1,\dots,z_{k-1}\}$.  We claim that  $y_k$ has no neighbours in $W$. Otherwise,  if $y_k$ has a neighbour, say $z$, in $W\backslash W_1$  then   $\bigcup _{1\leq i\leq k-1} Q_i$ together with a disjoint $P_3$ with center $y_k$ and two ends $u, z$  yield  $k\cdot P_3$, a contradiction.
If   $y_k$ has a neighbour  in $W_1$, without loss of generality, say
$x_1$, then    $\bigcup _{2\leq i\leq k-1} Q_i$ together with $zy_1z_1$ and $uy_kx_1$ will yield $k\cdot P_3$,  where $z$ is a neighbour of $y_1$ in $W\backslash W_1$.   It is  a contradiction. This implies that
 $N(y_k)\subseteq N(u)\cup \{u\}$ and $d(y_k)\leq k$.
By Claim,
\begin{eqnarray*}
  0 &\leq& \sum\limits_{v\in N(u)} d(v)- (k-1)d(u)-[(k-1)n-(k^2-k-1)] \\
  &=& \sum\limits_{i=1}^{k-1}d(y_i)+d(y_k)-(k-1)k-[(k-1)n-(k^2-k-1)]\\
   &\leq& (k-1)(n-1)+k-(k-1)k-[(k-1)n-(k^2-k-1)] \\
   &=& 0,
\end{eqnarray*}
implying that $d(y_1)=\cdots=d(y_k)=(n-1)$, $d(y_k)=k$, and $N(y_k)=\{u,y_1,\dots,y_k\}$.
Since $G$ is $k\cdot P_3$-free, $G-\{u,y_1,\dots,y_{k}\}$ is $P_3$-free. So  $G-\{u,y_1,\dots,y_{k}\}$
consists of independent edges and isolated vertices, which implies that $G\subseteq F_{n,k}$.
Then $\rho(G)\leq \rho(F_{n,k})$,  which implies that $\rho(G)=\rho(F_{n,k})$. By the Perron Fronbenius theorem and the extremality of $G$, we have  $G=F_{n,k}$.

\vspace{2mm}
{\bf Case 2:} $d(u)\geq3k-1$.
Then $G[N(u)]$ must be $(k-1)P_3$-free. Otherwise, $(k-1)\cdot P_3$ in $G[N(u)]$ and  a disjoint $P_3$ with center $u$ and two ends in $N(u)$ will yield  $k\cdot P_3$, a contradiction.  Similarly, $G-u$ is also $(k-1)\cdot P_3$-free.
By Theorem~\ref{edge extremal graph for kP3},
\[e(N(u))\leq\left\{
\begin{array}{llll}
 \vspace{1mm}
  \left(\begin{array}{c} 3k-4\\ 2\end{array}\right)+\Big\lfloor\frac{d(u)-3k+4}{2}\Big\rfloor, ~~~~~\mbox{for $3k-1\leq d(u)\leq5k-6$};\\
   \vspace{1mm}
  \left(\begin{array}{c} k-2\\ 2\end{array}\right)+(d(u)-k+2)(k-2)+\Big\lfloor\frac{d(u)-k+2}{2}\Big\rfloor, ~~~~\mbox{for $d(u)> 5k-6$};\\
\end{array}\right.
\]
and
$$e(G-u)\leq \left(\begin{array}{c} k-2\\ 2\end{array}\right)+(n-k+1)(k-2)+\bigg\lfloor\frac{n-k+1}{2}\bigg\rfloor.$$
If $3k-1\leq d(u)\leq5k-6$,
then
\begin{eqnarray*}
  && \sum\limits_{v\in N(u)} d(v) \leq d(u)+e(N(u))+e( G-u)\\
   &\leq&d(u)+\bigg( \left(\begin{array}{c} 3k-4\\ 2\end{array}\right)+\bigg\lfloor\frac{d(u)-3k+4}{2}\bigg\rfloor\bigg)+\\
   &&\bigg( \left(\begin{array}{c} k-2\\ 2\end{array}\right)+(n-k+1)(k-2)+
   \bigg\lfloor\frac{n-k+1}{2}\bigg\rfloor\bigg)\\
   &\leq&d(u) + \frac{d(u)+9k^2-30k+24}{2} +\frac{(2k-3)(n-1)-k^2+2k}{2}\\
    &=&(k-1)d(u)+(k-1)n-(k^2-k-1)-\frac{n+(2k-5)d(u)-10k^2+32k-25}{2} \\
    &\leq&(k-1)d(u)+(k-1)n-(k^2-k-1)-\frac{n+(2k-5)(3k-1)-10k^2+32k-25}{2} \\
    &=&(k-1)d(u)+(k-1)n-(k^2-k-1)-\frac{n-4k^2+15k-20}{2}\\
   &<& (k-1)d(u)+(k-1)n-(k^2-k-1),
\end{eqnarray*}
which contradicts to the claim.
If $d(u)> 5k-6$,  then
\begin{eqnarray*}
   &&\sum\limits_{v\in N(u)} d(v) \leq d(u)+e(N(u))+e( G-\{u\})\\
   &\leq&d(u)+\bigg(  \left(\begin{array}{c} k-2\\ 2\end{array}\right)+(d(u)-k+2)(k-2)+\bigg\lfloor\frac{d(u)-k+2}{2}\bigg\rfloor\bigg)+\\
   &&\bigg( \left(\begin{array}{c} k-2\\ 2\end{array}\right) +(n-k+1)(k-2)+\bigg\lfloor\frac{n-k+1}{2}\bigg\rfloor\bigg)\\
   &\leq &  d(u)+\frac{(2k-3)d(u)-k^2+2k}{2}+ \frac{(2k-3)(n-1)-k^2+2k}{2}\\
    &=&(k-1)d(u)+(k-1)n-(k^2-k-1)-\frac{n-d(u)+2k-1}{2} \\
    &\leq&(k-1)d(u)+(k-1)n-(k^2-k-1)-k\\
   &<& (k-1)d(u)+(k-1)n-(k^2-k-1),
\end{eqnarray*}
which also contradicts to the claim.
This completes the proof.
\end{Proof}

Now we are ready to prove Theorem~\ref{spectral radius for all linear forest}.

\vspace{2mm}
\begin{Proof} Theorem \ref{spectral radius for all linear forest} directly follows from
Lemmas~\ref{spectral radius for linear forest1}, \ref{spectral radius for all linear forest2} and \ref{spectral radius for kP3}.
\end{Proof}

\section{Proofs of Theorems~\ref{edge extremal bipartite graph}, \ref{extremal bipartite graph for kP3}, and Corollary~\ref{least eigenvalue}}
In this section, we will prove Theorems~\ref{edge extremal bipartite graph}, \ref{extremal bipartite graph for kP3}, and
 Corollary~\ref{least eigenvalue}.

\vspace{2mm}
\noindent{\bf Proof of Theorem~\ref{edge extremal bipartite graph}:}
Let $G$ be a $k\cdot P_3$-free bipartite graph of order $n\geq 11k-4$ with the maximum number of edges. Since  $K_{k-1,n-k+1}$ is $k\cdot P_3$-free,
we have $$e(G)\geq e(K_{k-1,n-k+1})=(k-1)(n-k+1).$$ We will prove the assertion by induction on $k$.

If $k=2$, then $e(G)\geq n-1$. Then $G$ contains $P_3$ as a subgraph. Otherwise $G$ consists of independent edges and isolated vertices and so $e(G)\leq \frac{n}{2}<n-1$, a contradiction.
 Clearly, every connected bipartite graph of order at least 3 contains $P_3$ as a subgraph. Since $G$ is $2\cdot P_3$-free, $G$ has exactly one connected component $H$ of order $p\geq 3$ and any of the remaining components (if any exists) is either an edge or an isolated vertex. Then $H$ is $2\cdot P_3$-free  and   $$e(H)\geq e(G)-\Big\lfloor\frac{n-p}{2}\Big\rfloor\geq n-1-\frac{n-p}{2}=\frac{n+p-2}{2}.$$
If $3\leq p\leq 4$, then $e(H)\geq 10$. If $p\geq5$, then $e(H)\geq11$.
Let $Q=v_1v_2\cdots v_l$ be the longest path of order $l$ in $H$. Since $P_3\subseteq H$ and  $H$ is $2\cdot P_3$-free, we have $3\leq l\leq 5$.
If $l=3$, then $H$ is a star $T_{p,0}$. On the other hand, $e(G)\geq n-1$. Hence $G=H=T_{n,0}$. If $l=4$, then  $V(H)\backslash V(Q)$ is an independent set and all vertices in $V(H)\backslash V(Q)$ are adjacent  to precisely one of $v_2$ and $v_3$. Hence  $H= T_{p,1}$.  On the other hand, $e(G)\geq n-1$. Then $G=H=T_{n,1}$.
If $l=5$, then $H-V(Q)$ consists of independent edges and isolated vertices, say $u_1v_1,\dots,u_qv_q,w_1,\dots,w_r$,
where $q+r\geq 1$, in which each vertex has at most one neighbor, which is $v_3$, in $Q$.
Since $H$  is a bipartite graph without containing  $2\cdot P_3$, it is easy to see that $H=T_{p,s}$, where $2\leq s\leq \lfloor\frac{p-1}{2}\rfloor$. On the other hand,  $e(G)\geq n-1$. Then $G=H=T_{n,s}$,
 $2\leq s\leq \lfloor\frac{n-1}{2}\rfloor$.
This completes the proof for $k=2$.

Suppose that  the assertion  holds for  $k-1\geq2$. Since
 $$e(G)\geq (k-1)(n-k+1)>(k-2)(n-k+2),$$
 By the induction hypothesis, we have $(k-1)\cdot P_3\subseteq G$.
 We have the following claim.

\vspace{2mm}
 {\bf Claim 1:} There exist $k-1$ vertices in $G$ with degree at least $3k-1$.

 In fact, take each $P_3=xyz$ in  $(k-1)\cdot P_3$. Then $G-V(P_3)$ must be $(k-1)P_3$-free, since $G$ is  $k\cdot P_3-$free.  By the induction hypothesis, $e(G-V(P_3))\leq(k-2)(n-k+2)$. Moreover,  $e(G[P_3])=2$ since $G$ is bipartite. Hence
\begin{eqnarray*}
  e(V(P_3),V(G)\backslash V(P_3)) &=& e(G)-e(G-V(P_3))-e(G[P_3]) \\
   &\geq &  (k-1)(n-k+1)-(k-2)(n-k+2)-2 \\
   &=& n-2k+1.
\end{eqnarray*}
Then  there exists a vertex in $P_3=xyz$ with degree at least $\frac{n-2k+1}{3}\geq 3k-1$.
  Therefore, for $(k-1)\cdot P_3$, there exist  $k-1$ vertices with degree at least $3k-1$. This finishes the proof of Claim 1.

   Let $U$ be a set of  $k-1$ vertices  in $G$ with degree at least $3k-1$.
   Then $G-U$ is $P_3$ free, i.e.,  $G-U$ consists of independent edges and isolated vertices.
Otherwise the $P_3$ in $G-U$ and  other $(k-1)\cdot P_3$ with all centers in $U$ and all ends in $V(G-U)\backslash V(P_3)$ will yield    $k\cdot P_3$ in $G$  since each vertex in $U$  has degree at least $3k-1$.
Further we have the following claim.

\vspace{2mm}
 {\bf Claim 2:} $U$ is an independent set.

 Suppose that $U$ is not an independent set.   Since $G$ is bipartite,   $G[U]$ has a bipartition $U=U_1\cup U_2$ with $e(U_1,U_2)\geq 1$.
   On the other hand,  $G-U$ has also a bipartition $V(G-U)=W_1\cup W_2$ with $|W_1|\leq |W_2|$. Without loss of generality, $G$ has a  bipartition $V(G)=(U_1\cup W_1)\cup (U_2\cup W_2)$.  Then
\begin{eqnarray*}
  e(G) &=& e(U_1,U_2)+e(W_1,W_2)+e(U_1,W_2)+e(U_2,W_1) \\
   &\leq & |U_1||U_2|+|W_1| +|U_1||W_2|+|U_2||W_1|\\
   &<&  |W_2||U_2|+|U_1||W_1|+|U_1||W_2|+|U_2||W_1|\\
   &=& (|U_1|+|U_2|)(|W_1|+|W_2)\\
   &=& (k-1)(n-k+1),
\end{eqnarray*}
where the first inequality holds because $G-U$ consists of independent edges and isolated vertices and
$|W_1|\leq |W_2|$, the second inequality holds because $1\leq|U_1|\leq k-2 \leq \frac{n-k+1}{2}\leq |W_2|$.  This contradicts to $e(G)\ge (k-1)(n-k+1)$. Hence Claim 2 holds.

In addition, we have the following Claim 3.

\vspace{2mm}

{\bf Claim 3:}  $V(G-U)$ is an independent set.

Since $G-U$ is $P_3$-free, $G-U$ consists of independent edges and isolated vertices.
Suppose that there is an edge $uv$ in $G-U$.  Since $G$ is a bipartite graph, $u$ and $v$ has no common neighbours in $U$.
Then $d(u)+d(v)\leq k-1$.
Then
\begin{eqnarray*}
e(G) &\le & e(\{u, v\})+e(\{u, v\}, V(G)-\{u, v\})+e(G-\{u, v\})
 \\
 &\le & 1+k-1+(k-1)(n-k+1-2)\\
& <&(k-1)(n-k+1).
 \end{eqnarray*}
 It is a contradiction.  So Claim 3 holds.  Then $G$ is a bipartite graph with bipartite parts $U$ and $V(G)-U$.  Moreover, $|U|=k-1$  and $e(G)\ge  (k-1)(n-k+1)$. Therefore $G=K_{k-1, n-k+1}$.
\QEDB

\vspace{4mm}

\noindent{\bf Proof of Theorem~\ref{extremal bipartite graph for kP3}:}
Let $G$ be a $k\cdot P_3$-free bipartite graph of order $n$ with maximum  spectral radius. Since $K_{k-1,n-k+1}$ is  $k\cdot P_3$-free, we have $$\rho(G)\geq \rho( K_{k-1,n-k+1})=\sqrt{(k-1)(n-k+1)}.$$
By Theorems~\ref{upper bound for bipartite graph} and \ref{edge extremal bipartite graph},
$$\rho(G)\leq \sqrt{e(G)}\leq \sqrt{(k-1)(n-k+1)}.$$
Then $$\rho(G)=\sqrt{e(G)}= \sqrt{(k-1)(n-k+1)}.$$
By  Theorem~\ref{edge extremal bipartite graph}, $G=K_{k-1,n-k+1}$.
\QEDB

\vspace{4mm}
\noindent{\bf Proof of Corollary~\ref{least eigenvalue}:}
By a result of Favaron et al. \cite{FMS}, $\lambda_n(G) \geq \lambda_n(H)$  for some spanning bipartite subgraph $H$. Moreover, the equality holds if and only if $G=H$, which can be deduced by its original proof.
 By Theorem~\ref{extremal bipartite graph for kP3},
$$\rho(H)\leq \sqrt{(k-1)(n-k+1)}$$ with equality if and only if $H=K_{k-1,n-k+1}$.
Since  the spectrum of a bipartite graph is symmetry \cite{LP},  $$\lambda_n(H)\geq-\sqrt{(k-1)(n-k+1)}$$ with equality if and only if $H=K_{k-1,n-k+1}$.  Thus
 we have $$ \lambda_n(G)\geq -\sqrt{(k-1)(n-k+1)}$$
with equality if and only if $G=K_{k-1,n-k+1}$.
 \QEDB

\section{Discussion}
It is known that if $G$ is a $K_{r+1}$-free graph of order $n$, then $e(G)\le e(T_{n, r})$ and $\rho(G)\le \rho(T_{n, r})$ (for example, see \cite{Nikiforov2,Turan}), where $T_{n, r}$ is
 a complete $r$-partite graph of order $n$ with  partite sets of cardinalities $\lfloor n/r\rfloor$ or $\lceil n/r\rceil$. Further spectral Tur\'{a}n theorem ($\rho(G)\le \rho(T_{n, r})$) is a slight better than  Tur\'{a}n theorem ($e(G)\le e(T_{n, r})$). In fact, if $\rho(G)\le \rho(T_{n, r})$, then
 $$e(G)\le \bigg\lfloor\frac{n\rho(G)}{2}\bigg\rfloor\le \bigg\lfloor\frac{n\rho(T_{n,r})}{2} \bigg\rfloor\le e(T_{n, r}).$$
But in general, $e(G)\le e(T_{n, r})$ does not imply that $\rho(G)\le \rho(T_{n, r})$.

On the other hand, from Theorems~\ref{edge extremal graph for linear forest}, \ref{edge extremal graph for kP3}, and
\ref{spectral radius for all linear forest}, we see that if $G$ is an $F$-free graph of order $n$,  where $F$ is a linear forest, then the extremal graph  maximizing the number of edges is the same as the extremal graph  maximizing the spectral radius. Based on the above results, we may propose the following problem.

\begin{Problem} For a given graph $H$, let $G$ be an   $H$-free graph of order $n$. If the extremal graph maximizing the number of edges is the same as the extremal graph  maximizing the spectral radius, say $G^*$, then what is relations between $e(G)\le e(G^*)$ and $\rho(G)\le \rho(G^*)$?
\end{Problem}

 Clearly, if $G^*$ is $T_{n, r}$, then $\rho(G)\le \rho(G^*)$ implies that $e(G)\le e(G^*)$. But if $G^*$ is the extremal graph for a linear forest $F$, we give several examples to illustrate that they do not have ``implication" relation, while  if $G^*$ is the bipartite extremal graph for $k\cdot P_3$, they have ``implication" relation.

\vspace{2mm}

\noindent{\bf Example 1:}   Let $G=(K_r\vee (K_1\cup K_{l-r}))\cup \overline{K}_{n-l-1}$, where $ \left(\begin{array}{c} l\\ 2\end{array}\right)+r=hn-\frac{h^2+h}{2}$ and $0\leq r<l$. Then
$$e(G)=hn-\frac{h^2+h}{2}=e(S_{n,h}).$$ Since
$$hn-\frac{h^2+h}{2}= \left(\begin{array}{c} l\\ 2\end{array}\right)+r< \left(\begin{array}{c} l+1\\ 2\end{array}\right)<\frac{(l+1)^2}{2},$$ we have $$l>\sqrt{2hn-h^2-h}-1.$$
Thus\begin{eqnarray*}
      \rho(G) &\geq& \rho(K_l)=l-1 >\sqrt{2hn-h^2-h}-2\\
      &>&\frac{h-1+\sqrt{4hn-(3h^2+2h-1)}}{2}\\
       &=& \rho(S_{n,h}).
    \end{eqnarray*}
  Then $e(G)\leq e(S_{n,h})$  and  $\rho(G)>\rho(S_{n,h})$.
\vspace{2mm}

\noindent{\bf Example 2:}
Let $G$ be a $2h$-regular graph of large  order $n$. Obviously,
$$\rho(G)=2h<\frac{h-1+\sqrt{4hn-(3h^2+2h-1)}}{2}=\rho(S_{n,h})$$ and
$$e(G)=hn>e(S_{n,h}).$$
Then  $\rho(G)\leq\rho(S_{n,h})$  and  $e(G)> e(S_{n,h})$.

\vspace{2mm}

\noindent{\bf Example 3:}  Let $G=(K_r\vee (K_1\cup K_{l-r}))\cup \overline{K}_{n-l-1}$, where $ \left(\begin{array}{c} l\\ 2\end{array}\right)+r=\Big\lfloor\frac{(2k-1)n-k^2+1}{2}\Big\rfloor$ and $0\leq r<l$. Then
$$e(G)=\bigg\lfloor\frac{(2k-1)n-k^2+1}{2}\bigg\rfloor=e(F_{n,k}).$$ Since
$$\bigg\lfloor\frac{(2k-1)n-k^2+1}{2}\bigg\rfloor= \left(\begin{array}{c} l\\ 2\end{array}\right)+r< \left(\begin{array}{c} l+1\\ 2\end{array}\right)<\frac{(l+1)^2}{2},$$ we have $$l>\sqrt{(2k-1)n-k^2}-1.$$
Thus
 $$\rho(G)\geq l-1>\sqrt{(2k-1)n-k^2}-2>\frac{k-1+\sqrt{4(k-1)n-(3k^2-2k-5)}}{2}\geq \rho(F_{n,k}).$$

\noindent{\bf Example 4:}
Let $G$ be a $2k$-regular graph of  large order $n$. Obviously,
$$\rho(G)=2k<\frac{k-1+\sqrt{4(k-1)n-(3k^2-2k-1)}}{2}<\rho(F_{n,k})$$ and
$$e(G)=kn>\bigg\lfloor\frac{(2k-1)n-k^2+1}{2}\bigg\rfloor =e(F_{n,k}).$$

\noindent{\bf Proposition 5:}
Let $G$ be  a  bipartite graph $G$ of order $n$. If  $e(G)\leq e(K_{k-1,n-k+1})$, then
$\rho(G)\leq\rho(K_{k-1,n-k+1})$.

\vspace{2mm}

Indeed,
by Theorems~\ref{upper bound for bipartite graph} and \ref{edge extremal bipartite graph},
$$ \rho(G)\leq \sqrt{e(G)}\leq\sqrt{(k-1)(n-k+1)}=\rho(K_{k-1,n-k+1}).$$

However, $\rho(G)\leq\rho(K_{k-1,n-k+1})$ does not imply that
$e(G)\leq e(K_{k-1,n-k+1})$.  For example, let $G$ be a $2k$-regular bipartite graph of  large even order $n$. Obviously,
$$\rho(G)=2k<\sqrt{(k-1)(n-k+1)}<\rho(K_{k-1,n-k+1})$$ and
$$e(G)=kn>(k-1)(n-k+1) =e(F_{n,k}).$$

From above discussion, there is an interesting phenomenon for (spectral) Tur\'{a}n type problems.  Spectral Tur\'{a}n type results imply the corresponding Tur\'{a}n type results for some given forbidden graphs, while Tur\'{a}n type results imply the corresponding spectral Tur\'{a}n type results for other given forbidden graphs,. This phenomenon may be worth to further investigate.

\end{document}